\documentclass[twoside, 10pt]{amsart}
\usepackage{hyperref}
\hypersetup{
colorlinks=true,
urlcolor=blue,
citecolor=blue}
\usepackage{amsfonts,amssymb,amsmath,amsgen,amsopn,amsbsy,graphicx,epsfig}
\usepackage[utf8]{inputenc}\usepackage[english]{babel}

\newcommand{\example}[1][]{\noindent \emph{Example. }}
\newcommand{\note}[1][]{\noindent \emph{Note. }}
\newtheorem{lemma}{Lemma}

\newtheorem{defi}{Definition}

\newtheorem{theorem}{Theorem}
\newtheorem{remark}{Remark}

\makeatletter
\@namedef{subjclassname@2020}{%
  \textup{2020} Mathematics Subject Classification}
\makeatother

\title{On the coordinates of minimal vectors in a Minkowski-reduced basis.}

\author{\'A. G.Horv\'ath }

\address{\'A. G.Horv\'ath, Dept. of Algebra and Geometry, Budapest University of Technology,
Egry J\'ozsef u. 1., Budapest, Hungary, 1111, \\ORCID iD: https://orcid.org/0000-0003-2371-4818 }
\email{ghorvath@math.bme.hu}

\date{\today}


\begin{document}

\begin{abstract}
Finding the shortest vectors in a lattice is an NP-hard problem, so low-dimensional results also play an essential role in lattice reduction theory. Using Ryskov's result for the admissible centerings and Tammela's result for determining the Minkowski-reduced form, we prove that the absolute values of the coordinates of a minimal vector on a six-dimensional Minkowski-reduced basis are less than or equal to three. To sharpen P. Tammela's work, we combine some lattice geometry arguments with the aforementioned theoretical results.
\end{abstract}

\maketitle

\section{Introduction}
\label{Sec:1}

A lattice is a discrete $\mathbb{R}^n$ subgroup. Any lattice has a lattice basis, i.e. a set of linearly independent lattice vectors such that any lattice vector is a linear combination of the basis vectors with integer coefficients. Each base has the same number of elements; the number of elements is the dimension of the lattice. A reduced basis contains nearly orthogonal and preferably short vectors. There are many different concepts of reduction. The essential classical reductions are Hermite \cite{hermite}, Minkowski \cite{minkowski}, Korkine-Zolotarev \cite{korkine-zolotarev}, and Venkov \cite{venkov}. Comparing these methods from a theoretical point of view is more complex than it might seem at first glance. The long story that Minkowski and Hermite reductions agree in up to six dimensions and give definite results in higher dimensions are in \cite{ryskov(1972)}. For a more detailed history of reduction theory, see the book \cite{gruber} and the literature contained therein.

A crucial task in geometric number theory is finding the shortest vector in a lattice, called the minimal vector of the lattice. All classical reduction theories assume a priori the solvability of this problem. The problem of finding minimum vectors connects the reduction theory of positive definite quadratic forms, the theory of admissible centering of the minimal parallelepiped and the determination problem of the Dirichlet-Voronoi cell of the lattice. The geometric structure of the minimum vectors of the lattice is also crucial in other mathematical theories, such as discrete geometry, coding theory, and the root lattice theory used in modern differential geometry.

This article focuses on the coordinates of minimal vectors written on the Minkowski-reduced basis. We recall the definition of the Minkowski reduction and summarize the theory of admissible centerings and the results needed to prove Theorem \ref{thm:coordinates}.

\section{Minkowski-reduced basis}
\label{Sec:5}

In our article, the $n$-dimensional \emph{lattice} $L[e_1,\ldots,e_n]$ is the $n$-dimensional real vector space $\mathbb{R}^n=\mathrm{Lin}[e_1,\ldots , e_n]$ is the set of linear combinations of integer coefficients that we form from the linearly independent vector system $\{e_1,\ldots, e_n\}$. The lattice is a system of linearly independent lattice vectors $\{f_1,\ldots ,f_k\}$ \emph{primitive system},
if the relation $\mathrm{Lin}[f_1,\ldots, f_k]\cap L[e_1,\ldots,e_n]=L[f_1,\ldots,f_k]$ is satisfied. Each lattice base is a primitive system with $n$ elements. We say that $a_1$ is the \emph{minimal vector} of the lattice (or one of its shortest vectors) if its length is equal to the minimum value of the lengths of the vectors of the lattice. Every subset of a primitive system is also a primitive system, and minimal vectors are primitive systems with one element.

\begin{defi}\label{def:successiveminima}
Let $a_1$ be one of the shortest vectors of the lattice, let $a_2$ be one of the shortest vectors forming a linearly independent vector system with $a_1$, and similarly define the vectors $a_i$ for $i=3,\ldots, n$. The system $\{a_1,\ldots,a_n\}$ is a so-called \emph{system of vectors giving successive minima} (or briefly system of successive minima). The lengths of these vectors are the \emph{successive minima} of the lattice.
\end{defi}

In general, the system of successive minima is not form a basis in a lattice. Minkowski-reduced basis (or basis reduced by Minkowski) can be defined as follows:

\begin{defi}\label{def:minkowskireducedbasis}
Denote by $a_1$ a minimum vector of the lattice, let $a_2$ be one of the shortest vectors such that $a_1$ and $a_2$ are linearly independent vectors, so that they together form a two-element primitive system. Continuing the process we define similarly the other elements of the ordered basis $\{a_1,a_2, \cdots a_n\}$ which we call a \emph{basis reduced by Minkowski}.
\end{defi}

Alternatively, we may look at the lattice $L$ as the regular linear transformation of $Y=\mathbb{Z}^n$. Hence, $L=AY$ and the squared lengths of the lattice vectors can be obtained by the formula $Q(x)=|Ax|^2$ for $x\in Y$. In this case, $Q$ is a positive definite symmetric quadratic form. For a lattice, we can correspond an equivalence relation of positive forms, $Q$ and $Q'$ are equivalent to each other if there is a unimodular transformation $U$ holding the property $Q'(x)=Q(U(x))$ for all $x\in \mathbb{R}^n$. The equivalent forms correspond to distinct bases of the same lattice. A positive quadratic form $Q$ is called \emph{Minkowski-reduced} if $Q(u)\geq Q(e_i)$ holds for $1\leq i\leq n$ and each lattice point $u=(u_1,\ldots ,u_n)\in Y$ with $\mathrm{g.c.d.}(u_i,\ldots ,u_n)=1$. All lattices have a Minkowski-reduced basis, implying that each positive quadratic form is equivalent to a Minkowski-reduced form.

\section{The admissible centering of a lattice}
\label{Sec:3}

The lattice $L[e_1,\ldots,e_n]$ of dimension $n$ is a \emph{centering} of the lattice $L'[a_1,\ldots,a_n]$ (of the same dimension) if $L'[a_1,\ldots,a_n]\subset L[e_1,\ldots,e_n]$. We have to define first the concept of \emph{admissible centering} of the lattice $L'$ in the case where there is a basis of $L'$ of shortest vectors. (So in the lattice $L'$ there is a \emph{minimum basis}). A lattice $L'$ with a minimum basis \emph{ can be entered in an admissible way} by the lattice $L$, if $\min L'=\min L$. Since there is no minimum basis in all lattices, we can reformulate the definition of admissible centerings, concentrating the idea of successive minima.

\begin{defi}\label{def:admissiblecenterings}
Take $n$ independent lattice vectors $\{a_1,\ldots,a_n\}$ which lengths are the successive minima of the lattice $L'[e_1,\ldots,e_n]$. We say that $L'[e_1,\ldots,e_n]$ is an \emph{admissible centering of} $L[a_1,\ldots,a_n]$ if for arbitrary sequence $1\leq i_1<\ldots <i_k\leq n$ of indices in the lattice $\mathrm{Lin} [a_{i_1},\ldots, a_{i_k}]\cap L'[e_1,\ldots,e_n]$ the successive minima are $|a_{i_1}|\leq \ldots \leq |a_{i_k}|$.
\end{defi}

Every linearly independent set $\{a_1,\ldots,a_n\}$ of lattice vectors defines a lattice parallelepiped
$$
\Pi[a_1,\ldots,a_n]:=\{\sum\limits_1^n x_ia_i : 0\leq x_i<1 \mbox{for all} i \}.
$$
Let $L[a_1,\ldots,a_n]$ be a lattice with minimum basis $\{a_1,\ldots,a_n\}$ and  $L'[e_1,\ldots,e_n]$ one of its admissible centerings. Then there are finitely many points $\{s_1,\ldots s_j\}$ of $L'[e_1,\ldots,e_n]$ in the parallelepiped $\Pi[a_1,\ldots,a_n]$  with the following properties:

\begin{itemize}
\item The vector system $\{a_1,\ldots,a_n,s_1,\ldots s_j\}$ generates the lattice $L'[e_1,\ldots,e_n]$.
\item $\min L[a_1,\ldots,a_n]=\min L'[e_1,\ldots,e_n]$.
\end{itemize}

In this case, we also say that the parallelepiped $\Pi$ has an admissible centering. Since in the above case $\Pi$ is spent by minimum vectors of the lattice, we call it \emph{minimum parallelepiped}.

They also used the modified version of this notation for a lattice formed by successive minimum vectors of the lattice. Ryskov in \cite{ryskov(1976)} proved that for every admissible centering of a parallelepiped based on the successive minimum vectors of a lattice, there is an affinity which sends to this (general) admissible centering to an admissible centering of a parallelepiped spending by a system of minimum vectors of the corresponding lattice. This result says that in a fixed dimension, all combinatorial types of admissible centering arise as admissible centering of a minimum parallelepiped. Hence, the relative connection between the two lattices can always be described as an admissible centering of a lattice with a minimum basis. If a given dimension, the admissible centerings are classified, then we know the possible connections between the given lattice $L'[e_1,\ldots,e_n]$ and its sublattices $L[a_1,\ldots,a_n]$. The Table \ref{table:centerings} of Ryskov contains the most important geometric data of the admissible centerings of a minimum parallelepiped. The points $\{s_1,\ldots s_j\}$ defining the centering of the parallelepiped $\Pi[a_1,\ldots,a_n]$ have rational coordinates with respect to the basis $\{a_1,\ldots,a_n\}$. In fact, if a basis of $L'[e_1,\ldots,e_n]$ is $\{f_1,\ldots,f_n\}$, then the matrix $A$ of the linear transformation sending $f_i$ to $e_i$ has integer elements. The coordinates of the points $\{s_1,\ldots s_j\}$ are also integers with respect to the basis $\{f_1,\ldots,f_n\}$. Hence, the basis changes $f_i$ to $e_i$ changes the coordinates of $\{s_1,\ldots, s_j\}$ such rational coordinates whose denominators divide the determinant of $A$. Hence, the least common multiple $U$ of the denominators of the coordinates of the points $\{s_1,\ldots s_j\}$ divides the determinant $V$ of $A$. $V$ is called by \emph{the volume (or index) of the centering}. The \emph{relevant rows} of an admissible centering are the coordinates of the vectors $\{s_1,\ldots s_j\}$ with respect to the basis $\{a_1,\ldots,a_n\}$. Table \ref{table:centerings} contains the data of the admissible centerings up to dimension $6$.

\begin{table}
$$
\begin{array}{|c|c|c|cccccc|}
\hline
  \mbox{dimension} & U & V & \multicolumn{6}{c|}{\mbox{relevant rows}} \\
\hline
\hline
  2 & 1 & 1 & 0,&0& & & & \\
\hline
  3 & 1 & 1 & 0,&0,&0& & & \\
\hline
  4 & 1 & 1 & 0,&0,& 0,& 0 & & \\
\hline
  4 & 2 & 2 & {1}/{2},&{1}/{2},&{1}/{2},&{1}/{2}& & \\
\hline
  5 & 1 & 1 & 0,&0,&0,&0,&0 & \\
\hline
  5 & 2 & 2 & {1}/{2},&{1}/{2},&{1}/{2},&{1}/{2},&0 & \\
\hline
  5 & 2 & 2 & {1}/{2},&{1}/{2},&{1}/{2},&{1}/{2},&{1}/{2}&  \\
\hline
  6 & 1 & 1 & 0,&0,&0,&0,&0,&0\\
\hline
  6 & 2 & 2 & {1}/{2},&{1}/{2},&{1}/{2},&{1}/{2},&0,&0\\
\hline
  6 & 2 & 2 & {1}/{2},&{1}/{2},&{1}/{2},&{1}/{2},&{1}/{2},&0 \\
\hline
  6 & 2 & 2 & {1}/{2},&{1}/{2},&{1}/{2},&{1}/{2},&{1}/{2},&{1}/{2} \\
 \hline
  6 & 2 & 4 &  {1}/{2},&{1}/{2},&{1}/{2},&{1}/{2},&0,&0\\
   &  &  & {1}/{2},&{1}/{2},&0,&0,&{1}/{2},&{1}/{2}\\
   &  &  & 0,&0,&{1}/{2},&{1}/{2},&{1}/{2},&{1}/{2}\\
 \hline
  6 & 3 & 3 & {1}/{3},&{1}/{3},&{1}/{3},&{1}/{3},&{1}/{3},&{1}/{3}\\
\hline
\end{array}
$$
\caption{Relevant rows for $n\leq 6$.}
\label{table:centerings}
\end{table}

Note that Ryskov also gave a similar table of dimension seven in \cite{ryskov(1976)}, and of dimension eight by Zakharova and Novikova in \cite{zakharova(1980)}, but there is no such in higher dimensions.

\section{Characterization of the Minkowski-reduced forms by a finite system of inequalities.}
\label{Sec:4}

We turn now to the third from the following three questions of reduction theory:
\begin{itemize}
\item How can we find a reduced basis in a given lattice?
\item How can we find a minimum or short vector in a given lattice?
\item How can we find all minimum vectors in a given lattice if we know one of its reduced bases?
\end{itemize}
The first two questions were solved in low-dimensional cases. The reason is that Tammela in \cite{tammela(1973)} up to dimension six gave a chance to search Minkowski reduced basis. We recall his argument. Every basis $\{e_1,\ldots,e_n\}$ of a lattice determines a point in the $N=:\binom{n+1}{2}$-dimensional Euclidean space by the coefficients of its symmetric Gram matrix. The set of all such points is an open convex cone in $\mathbb{R}^N$ with vertex the origin. In this cone, the Minkowski-reduced forms set a sub-domain, denoted by $\mathcal{M}$. The domain $\mathcal{M}\cup\{0\}$ is a closed convex polyhedral cone with apex $0$. The explicit description of $\mathcal{M}$ for $n\leq 6$ is known. One can specify a concrete finite system of inequalities $Q(u)\geq Q(e_i)$ that determine $\mathcal{M}$.
It consists of the $n-1$ inequalities $Q(e_{i+1})\geq Q(e_i)$ and inequalities $Q(u)\geq Q(e_i)$ for which $u$ can be transformed into a column of the matrix in Table \ref{table:tammela}:

\begin{table}
$$
      \left(
        \begin{array}{ccccccccc}
          1 & 1 & 1 & 1 & 1 & 1 & 1 & 1 & 1 \\
          1 & 1 & 1 & 1 & 1 & 1 & 1 & 1 & 1 \\
          0 & 1 & 1 & 1 & 1 & 1 & 1 & 1 & 1 \\
          0 & 0 & 1 & 1 & 1 & 1 & 1 & 1 & 1 \\
          0 & 0 & 0 & 1 & 2 & 1 & 1 & 2 & 2 \\
          0 & 0 & 0 & 0 & 0 & 1 & 2 & 2 & 3 \\
        \end{array}
      \right)
$$
\caption{Checking inequalities of reduction when $n\leq 6$.}
\label{table:tammela}
\end{table}
by permuting the coordinates of $u$ and omitting the signees of these coordinates, while $\mathrm{g.c.d}(u_i,\ldots,u_n)=1$. This result was announced by Minkowski for $n\leq 6$ (only the columns with nonzero coordinates $\leq n$ are taken in dimension $n$) and proved by Minkowski for $n\leq 4$ (\cite{minkowski_1}), Ryskov (in \cite{ryskov(1973)}) and later Afflerbach (in \cite{afflerbach}) for $n=5$ and by Tammela for $n=6$. Slightly refined and extended the result to the case $n=7$ by Tammela in \cite{tammela(1977)}.

This result says that $\{e_1,\ldots,e_n\}$ forms a Minkowski-reduced basis in its lattice $L[e_1,\ldots,e_n]$ if for the above finite number of inequalities $Q(u)\geq Q(e_i)$ hold with the vectors permitted by the table. Conversely, if we find a vector $u$ that is allowed by the table and with which the inverse inequality is satisfied, the base under consideration is not Minkowski-reduced; the lattice vector $u$ is a better choice for the reduction than that which is an element of the base. Finding an algorithm for a minimum vector in a lattice is hard, so finding a Minkowski-reduced basis in a lattice is also not easy. First, Beyer, Roof and Williamson in \cite{beyer(1971)}, and Beyer in \cite{beyer(1972)}) did an algorithm for the determination of a Minkowski-reduced basis. In practice, this algorithm is used only $n\leq 6$. The Minkowski reduction theory can be extended to integral $n\times m$ matrices concerning the module $\mathcal{M}$, which is the module over the integers generated by the columns of the given $n\times m$ matrix $A$. Zassenhaus and Ford indicated that the number of operations needed to reduce the matrix decreases for fixed $m$ if $n$ tends to infinity. Donaldson established this fact in \cite{donaldson} using the probability method.
A little bit later, Afflerbach and Groethe presented a new algorithm which is more practicable for higher dimensions and requires less computation time (see in \cite{afflerbach(1985)}). The algorithm presented in their paper is usable to determine Minkowski-reduced lattice bases of pseudo-random number generators up to dimension 20. We suggest for the interested reader the paper of Agrell, Eriksson, Vardy and Zeger (see \cite{agrell(2002)}) to a good review of the closest lattice point and reduced bases problems.

We note that Tammela used his method to determine the Dirichlet-Voronoi cell of a lattice point. In \cite{tammela(1975)}, the author proved that for a Minkowski-reduced basis, we can find the integer representation of the relevant vectors of the Dirichlet-Voronoi cell of the lattice among those lattice vectors which are can be transformed to a column of the matrix of Table \ref{table:tammeladvcell}
\begin{table}
$$
      \left(
        \begin{array}{ccccccccc||ccccccccc||ccccc}
          1 & 1 & 1 & 1 & 1 & 1 & 1 & 1 & 1 & 1 & 1 & 1 & 1 & 1 & 1 & 1 & 1 & 1 & 1 & 1 & 1 & 1 & 1 \\
          1 & 1 & 1 & 1 & 1 & 1 & 1 & 1 & 1 & 0 & 1 & 1 & 1 & 1 & 1 & 1 & 1 & 1 & 1 & 1 & 1 & 1 & 2 \\
          0 & 1 & 1 & 1 & 1 & 1 & 1 & 1 & 1 & 0 & 1 & 1 & 1 & 1 & 1 & 1 & 2 & 1 & 1 & 1 & 1 & 2 & 2 \\
          0 & 0 & 1 & 1 & 1 & 1 & 1 & 1 & 1 & 0 & 2 & 2 & 1 & 2 & 2 & 2 & 2 & 1 & 2 & 2 & 2 & 2 & 2 \\
          0 & 0 & 0 & 1 & 2 & 1 & 1 & 2 & 2 & 0 & 0 & 2 & 1 & 2 & 2 & 2 & 2 & 2 & 3 & 3 & 3 & 3 & 3 \\
          0 & 0 & 0 & 0 & 0 & 1 & 2 & 2 & 3 & 0 & 0 & 0 & 3 & 2 & 3 & 4 & 3 & m & 0 & 3 & 4 & 4 & 3 \\
        \end{array}
      \right)
$$
\caption{Checking inequalities of the relevants of the Dirichlet-Voronoi cell.}
\label{table:tammeladvcell}
\end{table}
by permuting the coordinates and omitting the signees of these coordinates if we exclude the last column $(1,1,1,1,2,m)^T$ of the second part. The division of the matrix is because the first part is agreed by Table \ref{table:tammela}, so this is needed to find the reduced basis. The second part (with the mentioned column) contains the vectors that may be shorter than the longest basis vector included in the reduced basis, which is required to specify the vector. In the $n$-dimensional case, consider only those columns of the table in which the number of non-zero coordinates does not exceed $n$.

Since every vector of minimum length is a relevant vector of the Dirichlet-Voronoi cell of the origin, the coordinates of these vectors concerning a Minkowski-reduced basis are included in the columns of the table, apart from the sign and the permutations of the coordinates.

We note that the book by Achill Sch\"urmann \cite{schurmann} contains more recent results on 7-dimensional Minkowski-reduced forms, but our proof would lead to a highly complicated discussion in the case of 7-dimensional lattices. Therefore, in this paper, we only focus on the six-dimensional case.

\section{The theorem} \label{sec:theorem}

Let $n\leq 6$ and assume that the lattice $L[e_1,\ldots,e_6]$ is spanned by a Minkowski-reduced basis. Combining the results on admissible centerings with the results of reduction theory, we can prove a sharpening of Tammela's result on the coordinates of the minimum vectors. Table \ref{table:tammeladvcell} shows that the maximum absolute value of the coordinates in $5$ or $6$ dimensions can be $3$ or $4$. In the following, we prove better, actually sharp boundaries.

\begin{theorem}\label{thm:coordinates}
Let $L$ be an $n$-dimensional lattice, where $n\leq 6$, $m=\sum\limits_{i=1}^nx_ie_i$ be an arbitrary minimal vector of $L$, where $\{e_1, \ldots,e_n\}$ is a Minkowski-reduced basis of the lattice. Then, for an arbitrary index $i$, the following inequalities are satisfied:
$$
|x_i|\leq
\left\{\begin{array}{ccl}
    1 & \mbox{if} & n=2,3 \\
   2 & \mbox{if} & n=4,5 \\
   3 & \mbox{if} & n=6.
 \end{array}
\right.
$$
\end{theorem}

First, we highlight some easy geometric observations which we continuously use in our proof. Since we prove the statement from lower to higher dimensions step by step, without loss of generality, we can assume that for the investigated minimum vector $m=\sum x_ie_i$ holds that $x_i\ne 0$ for all $i$. In the following, we assume that the volume of the basis-parallelepiped is $1$. Hence, the volume of the investigated parallelepipeds are integer.

\begin{lemma}\label{lem: B}
If a parallelepiped $\Pi[a_1,\ldots ,a_n]$ has a $k$-dimensional centered face of denominator $2$ then its volume is even.
\end{lemma}

\begin{proof} If the centered face is $\Pi[a_1,\ldots, a_k]$ then the coordinates of the vector $b=(1/2)(a_1+\ldots +a_k)\in L[e_1,\ldots,e_n]$ are integers, respectively. Hence the volume of $\Pi[b,a_2,\ldots,a_n]$ is also an integer. But
$\mathrm{vol}(\Pi[b,a_2,\ldots,a_n])=\det[b,a_2,\ldots,a_n]=(1/2)\det[a_1+\ldots +a_k,a_2,\ldots, a_n]=(1/2)\mathrm{vol}(\Pi[a_1,a_2,\ldots,a_n])$ proves the statement.
\end{proof}

\begin{lemma}\label{lem: C}
For all indices $i$ the inequality $|x_i|\leq 4$ holds.
\end{lemma}

\begin{proof} The minimum vectors of the lattice $L[e_1,\ldots,e_n]$ are relevant vectors of the Dirichlet-Voronoi cell of the origin, so by Theorem 2 in \cite{tammela(1975)} the absolute value of the coordinates of the minimum vectors with respect to a Minkowski-reduced basis can be found in that modification of Table \ref{table:tammeladvcell} in which the last column of the middle part is omitted. Hence for all $i$ we have $|x_i|\leq 4$.
\end{proof}

\begin{lemma}\label{lem: D}
Let $\{e_1,e_2, \ldots e_6\}$ be a Minkowski-reduced basis of the lattice then in the proof of the theorem we can assume that $|e_1|=|e_2|=|e_3|=|e_4|=|m|=1$.
\end{lemma}

\begin{proof}
We recall a result of G. Cs\'oka from his paper \cite{csoka(1978)}. It says that in a $n$-dimensional lattice with $n\leq 6$, the Minkowski-reduced basis is the lexicographical minima of the set of such bases, which are ordered by the increasing lengths of their elements. (Observe that this theorem is a consequence of the fact (is proved by Ryskov) that the domains of the Minkowski-reduced forms, Hermite-reduced forms, and $L^\star$-reduced forms are agreed for $n\leq 6$.) From this, the first four elements of a Minkowski-reduced basis give a primitive system of successive minimum vectors. So, Ryskov's affinity sends these vectors to a primitive system of minimum vectors of the image lattice without changing the coordinates of a minimum vector.
\end{proof}

\begin{lemma}\label{lem: E} Assume that the lattice $L$ is spanned by the Minkowski-reduced basis $\{e_1,\ldots,e_n\}$, where $n\leq 6$. Let $\{e_1,\ldots,e_{n-2},x\}$ be a primitive system with the vector $x=\sum\limits_{i=1}^n x_ie_i$, having the properties that $x_n \ne 0$. (Clearly, we have that $|e_1|\leq \ldots \leq |e_{n-1}|\leq |x|$.) Then there is a column $\alpha:=\alpha_1e_1+\ldots +\alpha_{n-1}e_{n-1}+ \alpha_ne_n$ in the $n$-dimensional part of Table \ref{table:tammela} for which also holds that $\alpha=\sum\limits_{i=1}^{n-2} y_ie_i+y_{n-1}x$, consequently simultaneously hold the equalities
$$
\alpha_{n-1}=y_{n-1}x_{n-1} \mbox{ and } \alpha_n=y_{n-1}x_{n}
$$
with such non-zero integers which absolute values are in Table \ref{table:tammela}.
\end{lemma}

\begin{proof}
With a suitable affinity $\varphi$ we define a lattice of form $L^\star=\varphi(L)$ in which $\{e_1\ldots e_{n-2},x\}$ is a primitive system (hence $L[e_1,\ldots e_{n-2},x]$ is a common complete sublattice both of $L$ and $L^\star$) and the minimum value of the vectors out of $L[e_1,\ldots e_{n-2},x]$ is greater than $|x|$. Then $\{e_1\ldots e_{n-2},\varphi(e_{n-1}\}$ in its lattice couldn't be a Minkowski-reduced basis ($\varphi(e_{n-1}\not \in L[e_1,\ldots e_{n-2},x]$ hence there is a column $\alpha:=\alpha_1e_1+\ldots +\alpha_{n-1}\varphi(e_{n-1})+\alpha_n \varphi(e_{n})$ in the $n$-dimensional part of Table \ref{table:tammela}, which exclude $\varphi(e_{n-1})$ from the possible Minkowski-reduced basis elements. But the vector in $L^\star$ corresponding to $\alpha $ has to be shorter then $\varphi(e_{n-1})$. Hence it is in the common layer of $L$ and $L^\star$ implying that $\alpha =y_1e_1+\ldots +y_{n-2}e_{n-2}+y_{n-1}x$ where $y_i$'s are all non-zero integers as we stated.

Now we prove the existence of an $L^\star$ with the above properties. Since  $\{e_1,\ldots,e_{n-2},x\}$ is a primitive system we have $\mathrm{g.c.d}(x_{n-1},x_n)=1$. By our assumption $x_n\ne =0$. Let an orthogonal decomposition of $z\in\mathbb{R}^n$ is $z=z_1+z_z$ where $z_1\in \mathrm{Lin}[e_1,\ldots,e_{n-2},x]$. If $x_{n-1}=0$ then $L[e_1,\ldots e_{n-2},x]=L[e_1,\ldots e_{n-2},e_n]$ holds and $e_{n-1}$ each vector is in the form  $v=v_1e_1+\ldots v_{n-1}e_{n-1}+v_{n-1}x+ke_{n-1}$ with a $k\in \mathbb{Z}$. Let the affinity is $\varphi(z):=z_1+(1+\varepsilon)z_2$ with an arbitrary positive $\varepsilon$. Then for $k\ne 0$ $\varphi(v)=v_1e_1+\ldots v_{n-1}e_{n-1}+v_{n-1}x+k\varphi(e_{n-1})$ and for $v\in L\setminus L[e_1,\ldots,e_{n-2},x]$ the minimal value of $|\varphi(v)|=\varphi(e_{n-1})$ since the basis $\{e_1\ldots e_n\}$ is a Minkowski-reduced one. Clearly, we can choose $\varepsilon$ on such a way that the condition $|\varphi(e_{n-1})|>|x|$ also holds. In the other case, when $x_{n-1}\ne =0$ we have $e_{n-1}=\frac{1}{x_{n-1}} \left(x-x^\star-x_ne_n\right)\in L[e_1,\ldots,e_n]$, where $x^\star\in L[e_1,\ldots,e_{n-2}]$. Since the coordinates of $x-x^\star$ are integers (with respect to the basis $\{e_1,\ldots,e_n$) we get that $|x_{n-1}|=1$. Hence $x=x^\star \pm e_{n-1}+x_ne_n$ with a non-zero integer $x_n$. Hence $L[e_1,\ldots e_{n-2},x,e_{n-1}]$ is an sublattice of $L$ of index $|x_n|$. Apply the above affinity to this lattice with such $\varepsilon$ that the minimum value of $|\varphi(v)|$ for $v\in L\setminus L[e_1,\ldots,e_{n-2},x]$ will be grater than $|x|$. Since $\{e_1,\ldots, e_{n-2},e_{n-1}\}$ is a primitive system of $L$ then it is also a primitive system of $L[e_1,\ldots e_{n-2},x,e_{n-1}]$. Thus $\{e_1,\ldots, e_{n-2},\varphi(e_{n-1})\}$ is a primitive system of $\varphi(L[e_1,\ldots e_{n-2},x,e_{n-1}])$ but couldn't be a Minkowski-reduced basis in this lattice because of $|\varphi(e_{n-1})|>|x|$, as the required affinity exists in this case, too.
\end{proof}

\begin{remark}\label{rem:minbasiscase} The theorem statement is evident if $\{e_1,\ldots,e_n\}$ is a basis consisting of minimal vectors. In fact,
$$
\mathrm{vol}(\Pi[e_1,\ldots,e_{i-1},m,e_{i+1}\ldots e_n])=|x_i|\mathrm{vol. }(\Pi[e_1,\ldots,e_n]).
$$
Using this equality, $|x_i|$ can only be one of the allowed volumes in Table \ref{table:centerings}. However, this is greater than three only in one case, when $n=6$ and the volume of the parallelepiped $\Pi[e_1,\ldots,e_{i-1},m,e_{i+1}\ldots e_n]$ is 4. This case can only occur if the lower-dimensional sublattice $L'[e_1,\ldots,e_{i-1},e_{i+1}\ldots e_n]$ contains a four-dimensional central $L''$ cubic-lattice whose index is $2$ in the lattice $L'[e_1,\ldots,e_{i-1},e_{i+1}\ldots e_n]$. Then, however, the volume equality takes the following form:
$$
4=\mathrm{vol}(\Pi[e_1,\ldots,e_{i-1},m,e_{i+1}\ldots e_n]\geq |x_i|\mathrm{ind}(L''/ L'[e_1,\ldots,e_{i-1},e_{i+1}\ldots e_n])=2|x_i|,
$$
i.e. the inequality $|x_i|\leq 2$ also exists.
\end{remark}

\begin{proof}[Proof of Theorem \ref{thm:coordinates}]
Let $\{e_1,\ldots, e_n\}$ be a Minkowski-reduced basis and $m=\sum x_ie_i$ is a minimum vector in the lattice $L=[e_1,\ldots, e_n]$. Assume that the volume of the basic parallelepiped $\Pi[e_1,\ldots, e_n]$ equals 1. By Lemma \ref{lem: D}, we can assume that the coordinates $x_i$ are non-zero and the first four elements of the basis are minimum vectors. By Remark \ref{rem:minbasiscase}, the statement is trivial in dimensions $n\leq 4$; hence, we have to investigate only the cases $n=5$ and $n=6$, respectively.

\begin{center}
{\bf n=5}
\end{center}

We cannot examine the coordinates similarly, so we must consider separate cases in the proof.
\begin{itemize}
\item The case of the fifth coordinate $x_5$.
Consider the four-dimensional minimum parallelepiped $\Pi[m,e_1,e_2,e_3]$. Its edges have length 1, so the corresponding vector system is a system of successive minima. Then $e_4$ be a fifth successive minima, so the lattice $L[e_1,\ldots, e_5]$ is an admissible centering of the lattice $L[m,e_1,e_2,e_3,e_4]$. Hence $\mathrm{vol}L[m,e_1,e_2,e_3,e_4]=|x_5|\leq 2$ by Table \ref{table:centerings}.

\item The case of the fourth coordinate $x_4$.
If $\{e_1,e_2,e_3,m\}$ {\bf is a primitive system} then there is a Minkowski-reduced basis with first four element $\{e_1,e_2,e_3,m\}$. We use Lemma \ref{lem: E} with the choice $x=m$, so the coefficients $x_i$ are non-zero integers. Since by the $5$-dimensional part of Table \ref{table:tammela} for all $i$ $|\alpha_i|\leq 2$, and since $|y_4|\geq 1$ we get $|x_4|\leq 2$, too.

Assume now that $\{e_1,e_2,e_3,m\}$ {\bf is not a primitive system} in $L[e_1,\ldots,e_5]$. Then
$$
L[e_1,\ldots, e_5]\cap \mathrm{Lin(e_1,e_2,e_3,m)}
$$
is an admissible centering of $L[e_1,e_2,e_3,m]$. Hence the lattice $L[e_1,\ldots, e_5]\cap \mathrm{Lin(e_1,e_2,e_3,m)}$ is a centered cubic lattice of dimension four and volume $2$. Since $e_4$ is a fifth vector with the minimal length for which $\{e_1,e_2,e_3,e_4^\star,e_4\}$ should form a basis (there is no admissible centering with volume 4), we can conclude by Lemma \ref{lem: D} that $e_5$ is also a minimum vector. We can finish the proof of this case using the result of the Remark \ref{rem:minbasiscase}, so in this case, we also hold $|x_4|\leq 2$.

\item Clearly, to prove that $|x_i|\leq 2$ for $i=1,2,3$ the argument of $i=4$ can be applied again.
\end{itemize}

\begin{center}
{\bf n=6}
\end{center}

In what follows, we will use the fact that the statement is true for $n\leq 5$. Let $\{e_1,\ldots,e_6\}$ and $m=\sum_{i=1}^6 x_ie_i$ be a Minkowski-reduced basis and a minimal vector with non-zero integer coefficients.

We distinguish two cases containing this proof's first and second paragraphs.

\begin{enumerate}
\item {\bf If $\{e_1,\ldots,e_4,m\}$ is a primitive system} then $|e_5|=1$ and $\Pi[e_1,\ldots, e_5,m]$ is a minimum parallelepiped with volume $|x_6|$. By Lemma \ref{lem: C} $|x_6|\leq 4$ and we have to prove that $|x_6|\neq 4$. By Table \ref{table:centerings}, the only lattice centred by volume 4 is the 6-dimensional cubic lattice. This centering has only one metric realization (see \cite{ryskov(1976)} Theorem 6). The shortest vector of the lattice linearly independent of $\{e_1,\ldots,e_5\}$ and completes this system to a basis must also be a minimal vector. Therefore, we have $|e_6|=1$. We can again use the last argument of Remark \ref{rem:minbasiscase} to show that all coordinates of $m$ are less than or equal to 3, as already mentioned.

\item {\bf If $\{e_1,\ldots,e_4,m\}$ is not a primitive system} then the face $\Pi[e_1,\ldots,e_4,m]$ is a centered facet of the parallelepiped $\Pi[e_1,\ldots,e_4,m,e_5]$. By Lemma \ref{lem: B} $|x_6|=\mathrm{vol}(\Pi[e_1,\ldots,e_4,m,e_5])$ is an even number. (Note that in this case the columns of Table \ref{table:centerings} cannot be directly referenced because $|e_5|> 1$ can also exist.) According to \ref{lem: C} Lemma, the value of $|x_6|$ is 2 or 4.

\begin{itemize}

\item Suppose $|x_6|=2$. Then there is a vector $e_5^\star \in L[e_1,\ldots,e_4,m]$ for which $\{e_1,\ldots,e_4,e_5^\star\}$ is a primitive system. Based on Table \ref{table:centerings}, there are two essentially distinct options for generating the vector, either $e_5^\star=1/2(e_1+ \ldots + e_4+m)+L[e_1,\ldots ,e_4]$ or $e_5 ^ \star=1/2(e_1+\ldots + e_{i-1} + e_{i+1}+ \ldots + e_4+m)+L[e_1,\ldots ,e_{i-1},e_{ i +1},\ldots e_4]$.

\begin{itemize}
\item In the first case, $e_5^\star=\sum_{i=1}^4((x_i+k_i)/2)e_i+(x_5/2)e_5+e_6$, where the coordinates are integers. Applying Lemma \ref{lem:  E} with $x:=e_5^\star$ we have an integer vector $\alpha $ for which by the integer coefficient $y_{n-1}$ hold simultaneously the equalities $\alpha_{5}=y_{5}(x_{5}/2)$ and $\alpha_{6}=y_{5}$. Since $|\alpha_6|=|y_5|\leq 3$ and $|\alpha_5|=|y_5||x_5/2|\geq 3$ we have the possibilities: $|y_5|=3$ and $|x_5|\leq 2$; $|y_5|=2$ and $|x_5|\leq 3$; or $|y_5|=1$ and $|x_5|\leq 6$, respectively. By Lemma \ref{lem: C}, we have to exclude only the case when $|x_5|=4$ holds. Assume indirectly that $|x_5|=4$ and consider the system $\{e_1,e_2,e_3,m,e_6\}$. Based on the integer coordinates of $e_5^\star$, $x_4$ is odd, therefore $e_4$ is not a vector of the lattice $\mathrm{Lin}[e_1,e_2,e_3,m,e_6]\cap L[e_1,\ ldots ,e_6 ]$. We can apply Lemma \ref{lem: E}. It follows that there is an integer vector
    $$
    \beta=\beta_1e_1+\beta_2 e_2+\beta_3 e_3+\beta_6 e_6+\beta_4 e_4'+\beta_5 e_5'\in \mathrm{Lin}[e_1,e_2,e_3,m,e_6]\cap L[e_1,e_2,e_3,e_6,e_4',e_5'],
    $$
    with coordinates $|\beta_i|\leq 3$. Here with the representation $m=x_1e_1+x_2e_2+x_3e_3+x_4e_4'+x_5e_5'+x_6e_6$, $\beta$ has the form $y_1e_1+y_2e_2+y_3e_3+y_4m+y_6e_6$ implying that from $|\beta_5|=|y_4x_5|\leq 3$ follows three possibilities for $|y_4|=k/4$. They correspond to the values $k=1,2$ or $3$, respectively. ($k=0$ implies that $\beta\in L[e_1,e_2,e_3,e_6]$ therefore we cannot choose  $e_6$ to the sixth element of a Minkowski-reduced basis). But $y_4x_4$ is also an integer; hence, $x_4$ is even, which is a contradiction with the fact that we proved earlier that $x_4$ is odd. Finally, we got that in this case, $|x_5|\leq 3$ holds. But $|x_5|=2$ because it is non-zero and even. Investigate now the first four coordinates. By the above argument $|y_4|=k/2$ where $k=1,2$ or $3$. Since $x_4$ is odd then $k=2$, $|y_4|=1$ and $3\geq |\beta_4|=|x_4|$ proves the required inequality for the fourth coordinate $x_4$, too. The first four coordinates are equivalent to our arguments, so in this case, the equality $|x_i|\leq 3$ holds for all $i$.

\item Secondly, assume that $e_5^\star=\sum_{i=2}^4((x_i+k_i)/2)e_i+(x_5/2)e_5+e_6$, hence $e_5^\star$ gives a admissible centering of the minimal 4-dimensional parallelepiped $\Pi[e_2,e_3,e_4,m]$. This means that $e_5^\star$ is a minimal vector, and $\{e_1,e_2,e_3,e_4,e_5^\star\}$ is a primitive system of minimal vectors. So $e_5$ is a minimal vector, too. Since $\{e_1,e_2,e_3,e_4,e_5,e_5^\star\}$ is a basis, then $e_6$ is a minimal vector, and we can apply the result of Remark \ref{rem:minbasiscase}. Hence, the absolute values of the coordinates of $m$ are less or equal to $3$ in this subcase, too.

\end{itemize}

\item To end our discussion, assume that $\{e_1,\ldots,e_4,m\}$ is not a primitive system and $|x_6|=4$. Then the parallelepiped $\Pi[e_1,\ldots,e_4,m,e_5]$ has volume $4$.  Let $a_5$ be a sixth vector whose length is the sixth successive minimum of the lattice concerning the system $\{e_1,\ldots,e_4,m\}$ of minimum vectors. The volume of the parallelepiped $\Pi[e_1 \ldots,e_4,m,a_5]$ is less or equal to 4 because its admissible centering is the lattice $L[e_1 \ldots,e_6]$.
    \begin{itemize}
    \item If $\mathrm{vol}(\Pi[e_1,\ldots,e_4,m,a_5])=2$ and the vector $e_5^\star$ centers the facet $\Pi[e_1,\ldots,e_4,m]$ then we get that $\{e_1,\ldots,e_4,a_5,e_5^\star\}$ is a basis. If $|e_5^\star|\leq |a_5|$ then we get $|a_5|\geq |e_5^\star|\geq |e_5|$ implying that $|a_5|=|e_5|$. Now $\{e_1,\ldots,e_4,m,e_5\}$ a system of successive minimum vectors, and so the centering of $\Pi[e_1,\ldots,e_4,m,e_5]$ is the admissible centering of the cubic lattice with index $4$. But there is no minimum vector with coordinate $4$ in this lattice, giving a contradiction. So we have $|e_5^\star|> |a_5|$ and so $\{e_1,\ldots,e_4,a_5,e_5^\star\}$ is an ordered basis. Hence, again, $|a_5|=|e_5|$ leads to the same contradiction as above.
     \item If $\mathrm{vol}(\Pi[e_1,\ldots,e_4,m,a_5])=4$, then the centering is combinatorially agreed with the centering of the $6$-dimensional cube with index $4$. Hence, a $4$-dimensional centered face by index $2$ contains $a_5$ as an edge vector. Let $x$ be the shortest half diagonal of this face. It is shorter than the half-diagonal of a brick with the same edge lengths and longer than $a_5$.
         Thus $|a_5|^2\leq |x|^2\leq 1/4(3+|a_5|^2)$ so $|a_5|^2=1$ and the successive minimum vector is a minimum vector. The lattice $L[e_1,\ldots,e_4,m,a_5]$ is the $6$-dimensional cubic lattice and $L[e_1,\ldots,e_4,e_5,e_6]$ is its admissible centering. In this lattice, the Minkowski reduced basis is minimal, so we can apply Remark \ref{rem:minbasiscase} again.
     \end{itemize}

\end{itemize}
\end{enumerate}
\end{proof}

\end{document}